\documentclass[12pt]{amsart}
\usepackage{amsmath}
\usepackage{graphicx}
\usepackage{hyperref}
\usepackage[latin1]{inputenc}
\usepackage{mathtools}
\usepackage{amsfonts}
\usepackage{amssymb}
\usepackage[T1]{fontenc}
\usepackage{amsthm}
\usepackage{fullpage}
\usepackage{enumitem}
\usepackage{bbm}
\usepackage[usestackEOL]{stackengine}
\usepackage{multirow}
\newtheorem{theorem}{Theorem}[section]

\newtheorem{lemma}[theorem]{Lemma}

\newtheorem{proposition}{Proposition}[section]
\theoremstyle{definition}
\newtheorem{definition}{Definition}[section]

\theoremstyle{remark}

\numberwithin{equation}{section}

\title{Further Investigation on Cyclotomic mapping Permutation Polynomials over Finite Fields}


\author{Suman Mondal}
\address{Department of Mathematical Sciences, Tezpur University, Tezpur, Assam, 784028, India}
\email{mondalmondalsuman@gmail.com}

\begin{document}
	%\vspace{.3cm}
	\begin{abstract}
		We explore the connection between cyclotomic mapping permutation polynomials and permutation polynomials of the form $x^rf(x^{\frac{q-1}{l}})$ over finite fields. We present a new necessary and a new sufficient condition to verify permutation behavior of such polynomials over finite field. As its application, for particular values of $r$, we point out some permutation trinomials of the form $P(x)=2x^{r+8}+x^{r+4}+2x^r \in \mathbb{F}_{13}[x]$, and work on few classes of permutation binomials. 
	\end{abstract}
	
	\maketitle
	
	\section{Introduction}
    Let $p$ ba a prime, $m\in \mathbbm{N}$, and $q=p^m$. A polynomial is called a permutation polynomial (PP) over a finite field (FF) $\mathbbm{F}_q$ if it induces a bijective mapping from $\mathbbm{F}_q$ to itself. Going back to 19-th century, Hermite and later Dickson pioneered the study of permutation polynomials over finite fields, and in recent years, the study of permutation polynomials have increased because of their applications and involvements in public key cryptosystems (\cite{levine1987some},\cite{lidl1984permutation}), RC6 block ciphers (\cite{rivest2001permutation}), Tuscan-$k$ arrays (\cite{chu2002circular}), Costas arrays (\cite{golomb2002periodicity}), among many others. permutation polynomials are also used in coding theory, for instance, permutation codes in power communications (\cite{chu2004constructions}), and interleavers in Turbo codes (\cite{sun2005interleavers}) etc. In some of these applications, the study of permutation polynomials over finite fields has also been extended to the study of permutation polynomials over finite rings and other algebraic structures. 
    
    Several classes of permutation polynomials are explored based on their applications mainly in coding theory and cryptography. Throughout this paper, we focus on such a class of polynomials with at least one zero root over a finite field. We consider a polynomial $P(x) \in \mathbbm{F}_q[x]$ such that $P(0)=0$. In that case, $P(x)$ is of the form $P(x)=x^rf(x^s)$, where $0<r<q-1$ and $q-1=ls$ for some positive integer $l$ and $s$, and $f(x)$ is an arbitrary polynomial over $\mathbbm{F}_q$ of degree $e>0$. 
    
    As in \cite{wang2007cyclotomic}, we use the $r$-th order cyclotomic mappings $f^{r}_{A_{0},\;A_{1},\;A_{2},\cdots, \;A_{l-1}}$ of index $l$ and reveal a simple and very useful connection between the polynomials of the form $P(x)=x^rf(x^s)$ and the $r$-th order cyclotomic mapping polynomials $f^{r}_{A_{0},\;A_{1},\;A_{2},\cdots, \;A_{l-1}}(x)$. That is, $P(x)=x^rf(x^s)=f^{r}_{A_{0},\;A_{1},\;A_{2},\cdots, \;A_{l-1}}(x)$, where $A_i=f( \xi^i )$ for $0\leq i \leq {l-1}$ and \;$\xi$\; is a primitive $l$-th roots of unity (Lemma \ref{rela}). In (\cite{wang2007cyclotomic}), we use two necessary conditions to check the permutation behavior of any given polynomial of the form $P(x)=x^rf(x^s)$. Those conditions help to identify the polynomials which are not permutation polynomials. In Theorem (\ref{NNC}), we present a new necessary condition for $P(x)=x^rf(x^s)$ to be a permutation polynomial over $\mathbbm{F}_q$.  We know that properties of $A_i$'s are crucial while discussing the permutation behavior of such polynomials. In this case, we use the index of $A_i$' s in $ \mathbbm{F}_q$ for the necessary condition.

    Using this new necessary condition and the existing results, in Theorem (\ref{Suff}), we present a new sufficient condition for $P(x)=x^rf(x^s)$ to be a permutation polynomial over $\mathbbm{F}_q$. This condition also involves the index of $A_i$'s. We include these results in Theorem 1 of (\cite{wang2007cyclotomic}), and in Theorem (\ref{main}), we present the necessary and sufficient conditions for $P(x)=x^rf(x^s)=f^{r}_{A_{0},\;A_{1},\;A_{2},\cdots, \;A_{l-1}}(x)$ to be a permutation polynomial over $\mathbbm{F}_q$. For particular values of $r$, we also present some permutation trinomials of the form $P(x)=2x^{r+8}+x^{r+4}+2x^r \in \mathbb{F}_{13}[x]$.
    
    Finally, we explore a few classes of permutation binomials of the form $x^r(x^{es}+1)$ where $s,l,r,e$ are some related positive integers. As an application, we characterize $x^r(x^{es}+1)$ in terms of the new necessary and sufficient condition, and the index of $A_i$'s.
    \section{Cyclotomic mapping permutation polynomials}
   Let $\gamma$ be a primitive element of  $\mathbbm{F}_{q}$,\;$q-1=ls$ for some $l$, $s$ $\in$ $\mathbbm{Z^+}$ and $C_0$ be the collection of all $l$-th powers of $\gamma$. As $ c^q=c $,  $\forall  c \in \mathbbm{F}_q $(\cite{lidl1994introduction}), then $$C_0=\{ \gamma^{lj} : j=0, 1, 2, \cdots, s-1\}.$$
    Now trivially $ C_0 $ is a subgroup of the cyclic group $ ( \mathbbm{F}_q^\ast , \cdot )$, so the quotient group $\mathbbm{F}_{q}^\ast / C_0$ exists with respect to multiplication, with index $l$. The elements of $\mathbbm{F}_{q}^\ast / C_0$ are called the cyclotomic cosets $C_i$  and  are defined as
    \begin{center}
    {$C_i=\gamma^i C_0,\;\; \forall i=0, 1, 2, \cdots, l-1.$}
    \end{center}
    Let $x \in C_i$ for some $i \in \{ 0, 1, 2, \cdots, l-1 \} $, then $x$ is of the form $\gamma^{i+lj}$  where $j \in \{ 0, 1, 2, \cdots, l-1 \}$.
    For $r\in \mathbbm{Z^+}$ and any $ {A_0}, A_1, A_2, \cdots, A_{l-1} \in \mathbbm{F}_{q}$, we define \textit{$r$-th order cyclotomic mapping $f^{r}_{A_{0},\;A_{1},\;A_{2},\cdots, \;A_{l-1}}$ of index $l$} from $\mathbbm{F}_{q}$ to itself, as \\
   $$f^r_{{A_0}, A_1, A_2, \cdots , A_{l-1}}(x)
   =\begin{cases}
       0 &  \text{if}  \ x=0\\
       A_i x^r & \text{if} \ x \in C_i, i=0,1, \cdots , l-1.
   \end{cases}$$
    \noindent $f^{r}_{A_{0},\;A_{1},\;A_{2},\cdots, A_{l-1}}$ is called  the \textit{$r$-th order cyclotomic mapping of least index $l$} if $l$ be the least positive integer such that the mapping can be written as cyclotomic mapping. The polynomial $f^{r}_{A_{0},\;A_{1},\;A_{2},\cdots, \;A_{l-1}}(x)$ over $\mathbbm{F}_q$ of degree at most $q-1$ representing cyclotomic mapping $f^{r}_{A_{0},\;A_{1},\;A_{2},\cdots, \;A_{l-1}}$, is called an \textit{$r$-th order cyclotomic mapping polynomial}. In particular, if $r=1$, the polynomial obtained is known as \textit{cyclotomic mapping polynomial}.\\
    Let $\xi=\gamma^s$, then $\xi$ is a primitive $l$-th roots of unity. Now for $i=0, 1, 2, \cdots, l-1$; we define $A_i=f (  \xi^i )$ where $\xi$ is a primitive $l$-th roots of unity.
    \begin{lemma}\label{rela}
        For any $r\in \mathbbm{Z^+}$, \;$x^rf(x^s)=f^{r}_{A_{0},\;A_{1},\;A_{2},\cdots, \;A_{l-1}}(x)$ where $A_i=f (  \xi^i )$ for $0\leq i\leq l-1$ and $\xi$ is a primitive $l$-th roots of unity.
    \end{lemma}
    \begin{proof}
        For $x=0$, the equality holds trivially. For $x \in \mathbbm{F}_q^\ast$, let $x\in C_i$ for some $i \in \{ 0, 1, 2, \cdots, l-1 \} $. Then $x$ is of the form $\gamma^{i+lj}$ for some $j \in \{ 0, 1, 2, \cdots, l-1 \} $.\\
    Now, $x^rf(x^s)=x^rf(\gamma^{s\; (i+lj)})=x^rf(\gamma^{is}\gamma^{(ls)j})=x^rf(\gamma^{is})=x^rf(\xi^i)=x^rA_i$, for $0\leq i\leq l-1$.\\
    Hence, $x^rf(x^s)=f^{r}_{A_{0},\;A_{1},\;A_{2},\cdots, \;A_{l-1}}(x)$ where $A_i=f (  \xi^i )$ for $0\leq i\leq l-1$,  $\xi$ is a primitive $l$-th roots of unity.
    \end{proof}
    Suppose $P(x)=x^rf(x^s)=f^{r}_{A_{0},\;A_{1},\;A_{2},\cdots, \;A_{l-1}}(x)$ is a permutation polynomial over $\mathbbm{F}_q$, then from \cite{wang2007cyclotomic}, we have $(r,s)=1$ and $A_i=f(\xi^i)\neq 0$, $\forall i=0, 1, 2, \cdots, l-1$. We know that these necessary conditions help to point out the polynomials showing no permutation behavior. Here we present such a necessary condition that also points out the polynomials with no permutation behavior.
    \begin{lemma}\label{FF} (\cite{lidl1994introduction})
        Let $n$ be a positive integer and $K$ be a field of characteristic $p$ where $p$ is a prime.
        If $p$ does not divide $n$, then $E^{(n)}$ is a cyclic group of order $n$ with respect to multiplication in $K^{(n)}$.
    \end{lemma}
    Here $K^{(n)}$ is the splitting field of $x^n-1$ over the field $K$ and $E^{(n)}$ is the collection of all $n$-th roots of unity  over $K$. Let $K=\mathbbm{F}_q$, then $E^{(l)}= \langle \xi \rangle= \{ 1, \xi, \xi^2, \cdots, \xi^{l-1} \}$ as $\xi$ is a primitive $l$-th roots of unity and $p$ does not divide $l$.
    \begin{definition}\label{def}
        Let $\gamma$ be a primitive element of $\mathbbm{F}_q$, then for any non-zero element $a$ in $\mathbbm{F}_q$, $a$ can be presented as $a=\gamma^b$ for some non-negative integer $b$. Index of $a$ in $\mathbbm{F}_q$ is denoted by ${Ind_\gamma}(a)$, and defined as 
    \begin{center}
    {${Ind_\gamma}(a)\;\equiv\;  b\; ( \bmod\; q - 1$).}
    \end{center}
    That is, ${Ind_\gamma}(a)$ is the residue class $b$ mod $q-1$ such that $a=\gamma^b$.\\
    For example, ${Ind_\gamma}(1)\;\equiv\;  0\; ( \bmod\; q - 1$).
    \end{definition}
    \begin{lemma}\label{index}
        Let $a, b \in \mathbbm{F}_q^\ast $ with $a\neq b$, then 
     \begin{enumerate}[label=(\roman*)]
 \item ${Ind_\gamma}(ab)\;\equiv\; {Ind_\gamma}(a)+{Ind_\gamma}(b) \; ( \bmod\; q - 1)$;
         \item ${Ind_\gamma}(a/b)\;\equiv\; {Ind_\gamma}(a)-{Ind_\gamma}(b) \; ( \bmod\; q - 1)$;
          \item ${Ind_\gamma}(a^{-1})\;\equiv\;  -\; {Ind_\gamma}(a)\; ( \bmod\; q - 1)$;
          \item ${Ind_\gamma}(a_1 a_2 \cdots a_k)\;\equiv\;  \sum\limits_{i=1}^{{k}}{Ind_{\gamma}(a_k)} \; ( \bmod\; q - 1)$ for $a_1,a_2, \cdots, a_k \in \mathbbm{F}_q^\ast$;
        \item ${Ind_\gamma}(a^k)\;\equiv\;  k {Ind_{\gamma}(a)} \; ( \bmod\; q - 1)$. 
     \end{enumerate}
    \end{lemma}
    Using Definition \ref{def}, Lemma (\ref{index}) can be proved trivially.
    \begin{theorem}\label{NNC}
        Suppose $q-1=ls$ where $l$, $s$ are positive integers, and $r \in \mathbbm{N}$. If $P(x)$$=x^rf(x^s) \in \mathbbm{F}_q[x]$ is a permutation polynomial, then $l \mid 2 Ind_\gamma(A_0A_1 \cdots A_{l-1})$. 
    \end{theorem}
    \begin{proof}
        Let $B={A_0}^s{A_1}^s\xi^r{A_2}^s\xi^{2r}\cdots A_{l-1}^s\xi^{r(l-1)}$. Then from Lemma (\ref{index}), we have
        \begin{align*}
         Ind_{\gamma}(B)\equiv\sum\limits_{i=0}^{l-1}{Ind_{\gamma}(A_i^s)+\{srl(l-1)/2\}} \; (\bmod\; q-1).   
        \end{align*}
        Suppose $P(x)=x^rf(x^s)=f^{r}_{A_{0},\;A_{1},\;A_{2},\cdots, \;A_{l-1}}(x)$ is a permutation polynomial over $\mathbbm{F}_q$. Then from Theorem 1 in \cite{wang2007cyclotomic}, $\mu_l=\{ {A_0}^s,{A_1}^s\xi^r,{A_2}^s\xi^{2r},\cdots, A_{l-1}^s\xi^{r(l-1)} \}$ is the set of all the distinct $l$-th root of unity.\\
        As $\mathbbm{F}_q$ is of characteristic $p$ and $p$ does not divide $l$, from Lemma (\ref{FF}), we have $E^{(l)}=\mu_l.$ That is $\{ {A_0}^s,{A_1}^s\xi^r,{A_2}^s\xi^{2r},\cdots, A_{l-1}^s\xi^{r(l-1)} \} = \{ 1, \xi, \xi^2, \cdots, \xi^{l-1} \}$.\\
        So, $B=\xi^{ \{l(l-1)/2\}}$. If $l$ is even, then $2Ind_{\gamma}(B)\equiv {0} \; (\bmod\; q-1)$. If $l$ is odd, then $Ind_{\gamma}(B)\equiv {0} \; (\bmod\; q-1)$.\\
        For $l \in \mathbbm{Z}^+$, we have $2Ind_{\gamma}(B)\equiv {0} \; (\bmod\; q-1)$. That is,
        \begin{align*}
           2Ind_{\gamma}(B)\equiv 2\sum\limits_{i=0}^{l-1}{Ind_{\gamma}(A_i^s)+\{srl(l-1)\}} \; (\bmod\; q-1). 
        \end{align*}
        So, $q-1=ls \mid 2sInd_\gamma(A_0A_1 \cdots A_{l-1})$. That is, $l \mid 2Ind_\gamma(A_0A_1 \cdots A_{l-1})$.    
    \end{proof}
    In case $l\geq 3$, Theorem (\ref{NNC}) is useful to point out that for some given $P(x)=x^rf(x^s)=f^{r}_{A_{0},\;A_{1},\;A_{2},\cdots, \;A_{l-1}}(x)$ is not a permutation polynomial over $\mathbbm{F}_q$. Consider $P(x)=x^3+x^7=x^3f(x^4) \in \mathbbm{F}_{13}[x]$ where $f(x)=x+1$. Then $r=3,q-1=12,s=4,l=3,\gamma=2,\xi=3$ with $A_0=2, A_1=4, A_2=10.$\\
    Now $2Ind_{\gamma}(A_0A_1A_2)\equiv {2} \; (\bmod\; 12)$, so $3 \nmid 2Ind_{\gamma}(A_0A_1A_2)$. Using Theorem (\ref{NNC}), we find that $P(x)=x^3+x^7$ is not a permutation polynomial over $\mathbbm{F}_{13}$.
    
    Below, we list the necessary conditions for $P(x)=x^rf(x^s)=f^{r}_{A_{0},\;A_{1},\;A_{2},\cdots, \;A_{l-1}}(x)$ to be a permutation polynomial over $\mathbbm{F}_q$.
    \begin{theorem}\label{Newnecess}
        Let $r(<q-1),s,l$ be positive integers such that $q-1=ls$. If $P(x)=x^rf(x^s)=f^{r}_{A_{0},\;A_{1},\;A_{2},\cdots, \;A_{l-1}}(x)$ is a permutation polynomial over $\mathbbm{F}_q$, then we have the following.
        \begin{enumerate}[label=(\roman*)]
            \item $(r,s)=1$;
            \item $A_i=f(\xi^i)\neq 0$, $\forall i=0, 1, 2, \cdots, l-1$;
            \item $l \mid 2Ind_\gamma(A_0A_1 \cdots A_{l-1}).$
        \end{enumerate}        
    \end{theorem}
    The above conditions are not sufficient to verify that $P(x)$ is a permutation polynomial over $\mathbbm{F}_{q}$. Consider $P(x)=x^3+2x=x(x^2+2) \in \mathbbm{F}_5$. Then $r=1, s=2, l=2, \gamma=2, \xi=4$ with $A_0=3, A_1=1$. Here $2Ind_{\gamma}(A_0A_1=3)\equiv {2} \; (\bmod\; 4)$, so $l=2 \mid 2Ind_\gamma(A_0A_1=3)$.\\
    All the conditions of Theorem (\ref{Newnecess}) are satisfied in this case, however, observe that $P(2)=P(4)=2$. So, $P(x)=x^3+2x$ is not a permutation polynomial over $\mathbbm{F}_{5}$.
    \section{Further results involving index}\label{sec3}
    In the previous section, we discussed the necessary conditions to be a permutation polynomial over finite fields. In this section, we obtain a sufficient condition using the necessary condition discussed in Theorem (\ref{NNC}). We explore the application of these new necessary and sufficient condition and explore permutation behavior of few classes of polynomial. We also point out few permutation trinomials over $\mathbbm{F}_{13}.$
    
    The result below presents some strong conditions to inspect the permutation behavior of polynomials of the form $P(x)=x^rf(x^s)$ over the finite field $\mathbbm{F}_q$.
    \begin{theorem} \label{old}\label{MT}\cite{wang2007cyclotomic}
     Let $p$ be a prime , $q=p^m$ for $m \in \mathbbm{Z^+}$, $q-1=ls$ for some $l, s \in \mathbbm{Z^+}$, $\gamma$  be a primitive element of $\mathbbm{F}_q$,  $\xi = \gamma^s$  be a primitive $l$-th root of unity, and $P(x)=x^rf(x^s)=f^{r}_{A_{0},\;A_{1},\;A_{2},\cdots, \;A_{l-1}}(x)$ be a polynomial over $\mathbbm{F}_q$  with $(r, s)=1$ and $A_i \neq 0$, $\forall \;i=1, 2, \cdots, {l-1}$. Then the following are equivalent:  
    \begin{enumerate}[label=(\roman*)]
          \item $ P(x)=x^rf(x^s) \text{ is a permutation polynomial over  } \mathbbm{F}_q$.
          \item $f^{r}_{A_{0},\;A_{1},\;A_{2},\cdots, \;A_{l-1}}(x) \text{ is a permutation polynomial over }    
          \mathbbm{F}_q$.
          \item $ A_iC_{ir} \neq A_jC_{jr} \text{ for any } i,j  \text{ with } 0 \leq i < j \leq l-1$.
          \item $ {Ind_\gamma}(A_i/A_j)\;\not\equiv \;r(j-i)\; ( \bmod\; l) \; \text{ for any } i,j  \text{ with } 0 \leq i < j \leq l-1.$ 
          \item $ \{ {A_0}, A_1\gamma^r, A_2\gamma^{2r}, \cdots,  A_{l-1}\gamma^{(l-1)r}\} \text{ is a system of distinct representatives of } \mathbbm{F}_{q}^\ast / C_0.$
          \item $ \{ A_0^s, A_1^s\xi^r, \cdots,  A_{l-1}^s\xi^{(l-1)r}\}=\mu_l \text{ is the collection of all distinct } l \text{-th roots of unity}.$
          \item $ \sum\limits_{i=0}^{l-1} \xi^{cri}A_i^{cs}=0,  \;\;\forall \;c=1, 2, \cdots, {l-1}$.
    \end{enumerate}
    \end{theorem}
    Next, using the condition discussed in Theorem (\ref{NNC}), we explore a sufficient condition similar to Theorem (\ref{old}) $(iv)$ for $P(x)=x^rf(x^s)=f^{r}_{A_{0},\;A_{1},\;A_{2},\cdots, \;A_{l-1}}(x)$  to be a permutation polynomial over $\mathbbm{F}_q$.
    \begin{theorem}\label{Suff}
        Let $p$ be a prime , $q=p^m$ for $m \in \mathbbm{Z^+}$, $q-1=ls$ for some $l, s \in \mathbbm{Z^+}$, $\gamma$  be a primitive element of $\mathbbm{F}_q$,  $\xi = \gamma^s$  be a primitive $l$-th root of unity and $P(x)=x^rf(x^s)=f^{r}_{A_{0},\;A_{1},\;A_{2},\cdots, \;A_{l-1}}(x)$ be a polynomial over $\mathbbm{F}_q$  with $(r, s)=1$, $A_i \neq 0$ $\forall \;i=1, 2, \cdots, {l-1}$, and $l \mid 2Ind_\gamma(A_0A_1 \cdots A_{l-1})$. Then the following are equivalent.
        \begin{enumerate}[label=(\roman*)]
            \item $P(x)=x^rf(x^s)=f^{r}_{A_{0},\;A_{1},\;A_{2},\cdots, \;A_{l-1}}(x)$ is a permutation polynomial over $\mathbbm{F}_q$.
            \item $2{Ind_\gamma}(A_0A_1\cdots A_{i-1}A_{i+1} \cdots A_j^2 \cdots A_{l-1})\;\not\equiv \;2r(i-j)\; ( \bmod\; l)$ for any $i,j$ with $0\leq i<j\leq l-1$.
        \end{enumerate}
    \end{theorem}
    \begin{proof}
        Let $P(x)$ be a permutation polynomial over $\mathbbm{F}_q$. Then using Theorem (\ref{MT}), we have 
        \begin{align*}
            & {Ind_\gamma}(A_i/A_j)\;\not\equiv \;r(j-i)\; ( \bmod\; l) \; \text{ for any } i,j  \text{ with } 0 \leq i < j \leq l-1.
        \end{align*}
        Now given that $2Ind_{\gamma}(A_0A_1 \cdots A_{l-1})\equiv {0} \; (\bmod\; l)$. Using Lemma (\ref{index}), we have\\
          $2[{Ind_\gamma}(A_i/A_j)+{Ind_\gamma}(A_0A_1\cdots A_{i-1}A_{i+1} \cdots A_j^2 \cdots A_{l-1})]\;\equiv \;0\; ( \bmod\; l)$ for any $i,j$ with $0\leq i<j\leq l-1$.\\
          That is, $2{Ind_\gamma}(A_0A_1\cdots A_{i-1}A_{i+1} \cdots A_j^2 \cdots A_{l-1})\;\equiv \;-2{Ind_\gamma}(A_i/A_j)\; ( \bmod\; l)$.\\
          So, $2{Ind_\gamma}(A_0A_1\cdots A_{i-1}A_{i+1} \cdots A_j^2 \cdots A_{l-1})\;\not\equiv \;2r(i-j)\; ( \bmod\; l)$ for any $i,j$ with $0\leq i<j\leq l-1$.

          Conversely, let condition $(ii)$ be true.\\
          Suppose $P(x)$ is not a permutation polynomial over $\mathbbm{F}_q$. Then from Theorem (\ref{MT}), for some $i,j$ with $0\leq i<j \leq l-1,$ we have ${Ind_\gamma}(A_i/A_j)\;\equiv \;r(j-i)\; ( \bmod\; l).$\\
          As $2Ind_{\gamma}(A_0A_1 \cdots A_{l-1})\equiv {0} \; (\bmod\; l)$, we have\\ $2[{Ind_\gamma}(A_i/A_j)+{Ind_\gamma}(A_0A_1\cdots A_{i-1}A_{i+1} \cdots A_j^2 \cdots A_{l-1})]\;\equiv \;0\; ( \bmod\; l)$. That is,\\
          $2{Ind_\gamma}(A_0A_1\cdots A_{i-1}A_{i+1} \cdots A_j^2 \cdots A_{l-1})\;\equiv \;2r(i-j)\; ( \bmod\; l)$ for some $i,j$ with $0\leq i<j\leq l-1$, which is a contradiction.\\
          Hence, $P(x)$ is a permutation polynomial over $\mathbbm{F}_q$.   
    \end{proof}
    Using Theorem (\ref{NNC}) and Theorem (\ref{Suff}), below we refine Theorem (\ref{old}).
    \begin{theorem}\label{main}
        Let $p$ be a prime , $q=p^m$ for $m \in \mathbbm{Z^+}$, $q-1=ls$ for some $l, s \in \mathbbm{Z^+}$, $\gamma$  be a primitive element of $\mathbbm{F}_q$,  $\xi = \gamma^s$  be a primitive $l$-th root of unity and $P(x)=x^rf(x^s)=f^{r}_{A_{0},\;A_{1},\;A_{2},\cdots, \;A_{l-1}}(x)$ be a polynomial over $\mathbbm{F}_q$  with $(r, s)=1$, $A_i \neq 0$ $\forall \;i=1, 2, \cdots, {l-1}$, and $l \mid 2Ind_\gamma(A_0A_1 \cdots A_{l-1})$. Then the following are equivalent.
        \begin{enumerate}[label=(\roman*)]
          \item $ P(x)=x^rf(x^s) \text{ is a permutation polynomial over  } \mathbbm{F}_q$.
          \item $f^{r}_{A_{0},\;A_{1},\;A_{2},\cdots, \;A_{l-1}}(x) \text{ is a permutation polynomial over }    
          \mathbbm{F}_q$.
          \item $ A_iC_{ir} \neq A_jC_{jr} \text{ for any } i,j  \text{ with } 0 \leq i < j \leq l-1$.
          \item $ {Ind_\gamma}(A_i/A_j)\;\not\equiv \;r(j-i)\; ( \bmod\; l) \; \text{ for any } i,j  \text{ with } 0 \leq i < j \leq l-1.$
          \item $2{Ind_\gamma}(A_0A_1\cdots A_{i-1}A_{i+1} \cdots A_j^2 \cdots A_{l-1})\;\not\equiv \;2r(i-j)\; ( \bmod\; l)$ for any $i,j$ with $0\leq i<j\leq l-1$.
          \item $ \{ {A_0}, A_1\gamma^r, A_2\gamma^{2r}, \cdots,  A_{l-1}\gamma^{(l-1)r}\} \text{ is a system of distinct representatives of } \mathbbm{F}_{q}^\ast / C_0.$
          \item $ \{ A_0^s, A_1^s\xi^r, \cdots,  A_{l-1}^s\xi^{(l-1)r}\}=\mu_l \text{ is the collection of all distinct } l \text{-th roots of unity}.$
          \item $ \sum\limits_{i=0}^{l-1} \xi^{cri}A_i^{cs}=0,  \;\;\forall \;c=1, 2, \cdots, {l-1}$.
    \end{enumerate}
    \end{theorem}
    \begin{theorem}\label{lab}
        Let $p$ be a prime number, $q=p^m$ for some $m \in \mathbbm{Z^+}$, $q-1 = 3s$ for some $  s \in \mathbbm{Z^+}$. Assume $f(x) \;\equiv\; ax^2+bx+c\; ( \bmod\; x^3-1)$ such that $a^2+b^2+c^2-ab-bc-ca =1$. Then $P(x)=x^rf(x^s)=f^{r}_{A_{0},\;A_{1},\;A_{2}}$$(x)$ is a permutation polynomial over $\mathbbm{F}_q$ if and only if $(r,s)=1$, $A_0^s = 1$, $3\mid Ind_{\gamma}(A_0)$, $3 \nmid \{r+Ind_{\gamma}(A_2^2)\}$ where $\xi^3=1$ and $A_i=f (  \xi^i ) \neq 0$, $\forall$$i=0, 1, 2$.
    \end{theorem}
    \begin{proof}
        As $a^2+b^2+c^2-ab-bc-ca =1$, then trivially $A_1A_2=1$. If $P(x)$ is a permutation polynomial over $\mathbbm{F}_q$, then $\prod_{x \in \mathbbm{F}_q^\ast}$  $P(x) = -1$   implies $A_0^s=1$. So $A_i=f (  \xi^i ) \neq 0$, $\forall$$i=0, 1, 2$.\\
        From Theorem (\ref{main}), we have that $P(x)$ is a permutation polynomial over $\mathbbm{F}_q$ if and only if $(r,s)=1, 3\mid Ind_{\gamma}(A_0)$ and $\mu_3=\{ 1,A_1^s\xi^r, A_2^s\xi^{2r}\}$ is the collection of all distinct  $3$ -th roots of unity. We observe that every element of $\mu_3$ is a 3-th root of unity.\\
        From (\cite{wang2007cyclotomic}), we know that Showing $\mu_3$ is the collection of all distinct  $3$ -th roots of unity is equivalent with $A_1^s\xi^r \neq A_2^s\xi^{2r}$. Now 
        \begin{align*}
            & A_1^s\xi^r =A_2^s\xi^{2r}\\
 \Leftrightarrow & \; s{Ind_\gamma}(A_1/A_2)\;\equiv \;rs\; ( \bmod\; q-1)\\
 \Leftrightarrow & \; 2{Ind_\gamma}(A_0A_1A_2/A_0A_2^2)\;\equiv \;2r\; ( \bmod\; 3) \\
 \Leftrightarrow & \; 2s{Ind_\gamma}(A_0A_2^2)\;\equiv \;-2rs\; ( \bmod\; q-1)\\
 \Leftrightarrow & \; 2s{Ind_\gamma}(A_2^2)\;\equiv \;-2rs\; ( \bmod\; q-1)\\
 \Leftrightarrow & \; 2{Ind_\gamma}(A_2^2)\;\equiv \;-2r\; ( \bmod\; 3)\\
 \Leftrightarrow & \; 3 \mid r+Ind_{\gamma}(A_2^2).
        \end{align*}
        So, $A_1^s\xi^r \neq A_2^s\xi^{2r}$ is equivalent with $3 \nmid r+Ind_{\gamma}(A_2^2)$.\\
        Hence the theorem.
    \end{proof}
    \noindent\textbf{Example 3.1.} Consider $P(x)=2x^9+x^5+2x=x(2x^8+x^4+2)=xf(x^4) \in \mathbbm{F}_{13}[x]$, where $f(x)=2x^2+x+2$. Then $r=1, q-1=12, l=3, s=4, \gamma=2, \xi=3$ with $(r,s)=1$ and $A_0=f(1)=5, A_1=f(3)=10, A_2=f(9)=4, A_2^2=3, A_0^4=1.$\\
    Here ${Ind_2}(A_0)\;\equiv \;9\; ( \bmod\; 12)$ and ${Ind_2}(A_2^2)\;\equiv \;4\; ( \bmod\; 12)$, so $3\mid Ind_{\gamma}(A_0)$ and $3 \nmid \{r+Ind_{\gamma}(A_2^2)\}$.\\
    Using Theorem (\ref{lab}), $P(x)$ is a permutation polynomial over $\mathbbm{F}_{13}.$\\
    Again, ${Ind_2}(A_1^2A_2=10)\;\equiv \;10\; ( \bmod\; 12)$, ${Ind_2}(A_1A_2^2=4)\;\equiv \;2\; ( \bmod\; 12)$, ${Ind_2}(A_0^2A_2=2)\;\equiv \;1\; ( \bmod\; 12)$. So, $2{Ind_2}(A_1^2 A_2=10)\;\not\equiv \;\{2\cdot1 \cdot (0-1)\}\; ( \bmod\; 3)$, $2{Ind_2}(A_1A_2^2=4)\; \not\equiv \;\{2\cdot 1 \cdot (0-2)\}\; ( \bmod\; 3)$, $2{Ind_2}(A_0 A_2^2=2)\;\not\equiv \;\{2\cdot1 \cdot (0-1)\}\; ( \bmod\; 3)$.\\
     Using Theorem (\ref{Suff}), $P(x)$ is a permutation polynomial over $\mathbbm{F}_{13}.$
     \begin{proposition}
         For $r=1,3,7,9$; $P(x)=2x^{r+8}+x^{r+4}+2x^r$ is a permutation polynomial over $\mathbbm{F}_{13}$.
     \end{proposition}
     \begin{proof}
         Here $P(x)=2x^{r+8}+x^{r+4}+2x^r=x^r(2x^8+x^4+2)=x^rf(x^4) \in \mathbbm{F}_{13
         }[x]$, where $f(x)=2x^2+x+2$. Taking $q-1=12, l=3, s=4, \gamma=2, \xi=3$, we have $A_0=f(1)=5, A_1=f(3)=10, A_2=f(9)=4, A_2^2=3, A_0^4=1.$\\
    Here ${Ind_2}(A_0)\;\equiv \;9\; ( \bmod\; 12)$ and ${Ind_2}(A_2^2)\;\equiv \;4\; ( \bmod\; 12)$.\\
    From Theorem (\ref{lab}), $P(x)$ is a permutation polynomial over $\mathbbm{F}_q$ if and only if $(r,4)=1$ and $3 \nmid r+4$ where $0<r<12$.\\
    Hence, for $r=1,3,7,9$; $P(x)=2x^{r+8}+x^{r+4}+2x^r$ is a permutation polynomial over $\mathbbm{F}_{13}$.  
     \end{proof}
     \section{few classes of permutation binomials}
     In previous sections, we explored some necessary and sufficient conditions for $P(x)=x^rf(x^s)$ to be a permutation polynomial over $\mathbbm{F}_q$. As an application, we now focus on the polynomial of the form $P(x)=x^r(x^{es}+1) \in \mathbbm{F}_q[x]$ where $0<r<q-1, q-1=ls,$ and $e\in \mathbbm{N}$ with $(e,l)=1$. From (\cite{wang2007cyclotomic}), in this case we have $l$ is odd and $s$ is even. We consider $l\geq 3$. We also discuss the permutation behavior of a subclass of $P(x)$ over $\mathbbm{F}_q$. 
     \begin{theorem}\label{pbn}
       Let $p$ be an odd prime, and $q=p^m$ for $m \in \mathbbm{N}$. Assume $l,r,s,e \in \mathbbm{N}$ sohat $l(\geq 3)$ is odd, $(l,e)=1$, and $q-1=ls$. If $P(x)=x^r(x^{es}+1)$ is a permutation binomial over $\mathbbm{F}_q$ then $(r,s)=1$, $p \mid 2^s-1, l \nmid 2r+es$.       
     \end{theorem}
     \begin{proof}
         $(r,s)=1$ is trivial.
         
         As $l$ is odd and $(e,l)=1$, $\xi^e$ is also a primitive $l$-th root of unity and ${\displaystyle \prod_{i=0}^{l-1} \xi^{ei}}=1$.\\
         Now ${\displaystyle \prod_{i=0}^{l-1} A_i}={\displaystyle \prod_{i=0}^{l-1} (\xi^{ei}+1)}={\displaystyle \prod_{i=0}^{l-1} (1-(-\xi^{ei}))}=1-(-1)=2.$ That is, $A_1A_2 \cdots A_{l-1
         }=1$.\\
         From Theorem (\ref{Newnecess}) (iii), we have $l \mid Ind_\gamma(A_0A_1 \cdots A_{l-1}).$ So $l \mid Ind_\gamma(2)$, and for $\xi=\gamma^s$
         \begin{center}
            ${Ind_\gamma}(2^s)\;\equiv \;0\; ( \bmod\; q-1) $
         \end{center}
         Hence $ 2^s\;\equiv \;1\; ( \bmod\; p) $, that is, $p \mid 2^s-1.$

         Suppose $l \mid 2r+es$. As $l$ is odd and $s$ is even, we have $2l\mid 2r+es$.\\Now $l$ is odd and $l \mid q-1$. So we can find $\eta \in \mathbbm{F}_q^\ast$ such that $\eta^2=\xi$. By Theorem (\ref{main}) $(viii)$, we have $ \sum\limits_{i=0}^{l-1} \xi^{cri}A_i^{cs}=0,  \;\;\forall \;c=1, 2, \cdots, {l-1}$. That is, $\sum\limits_{i=0}^{l-1} \eta^{(2r+es)ci}(\eta^{ei}+\eta^{-ei})^{cs}=0,  \;\;\forall \;c=1, 2, \cdots, {l-1}$. So $\sum\limits_{i=0}^{l-1} (\eta^{ei}+\eta^{-ei})^{cs}=0,  \;\;\forall \;c=1, 2, \cdots, {l-1}$.\\
         As each $(\eta^{ei}+\eta^{-ei})^{s}$ is an $l$-th root of unity, using Lemma 2 in (\cite{wang2007cyclotomic}), $(\eta^{ei}+\eta^{-ei})^{s}$ are all distinct  for all $i=0,1, \cdots ,l-1$. However, as $s$ is even, we have $(\eta^{ei}+\eta^{-ei})^{s}=(\eta^{(l-i)e}+\eta^{-(l-i)e})^{s}$ which is a contradiction.\\
         Hence $l \nmid 2r+es$.         
     \end{proof}
     \begin{theorem}\label{PB}
         Let $p$ be an odd prime, and $q=p^m$ for $m \in \mathbbm{N}$. Assume $r,s,e \in \mathbbm{N}$ such that $(3,e)=1$ with $q-1=3s$. Then $P(x)=x^r(x^{es}+1)$ is a permutation polynomial over $\mathbbm{F}_q$ if and only if $(r,s)=1,  3 \mid 2^s-1, 3\nmid 2r+es,   3\nmid r+es,$ and $3\nmid r+2es$.  
     \end{theorem}
    \begin{proof}
        We have $P(x)=x^r(x^{es}+1)=x^rf(x^s)\in \mathbbm{F}_q[x]$ where $f(x)=x^e+1$. Suppose $\xi$ is a primitive $3$-th root of unity, and $\xi=\gamma^s$ with $A_i=f(\xi^i) \;\forall i=0,1,2$, then trivially $A_i\neq 0$.\\
        Now $A_0=2, A_1=\xi^e+1, A_2=\xi^{2e}+1$ with $A_1A_2=(\xi^e+1)(\xi^{2e}+1)=A_1+A_2$.\\
        As $\xi$ is a primitive $3$-th root of unity, we have $\xi^{2e}+\xi^e+1=0.$
        That is, $A_1^2=\xi^e.$ \\
        Now $A_1A_2=A_1+A_2$ implies $A_1(A_2-1)=A_2$ and $A_2(A_1-1)=A_1$. That is, $A_1\xi^{2e}=A_2$, $A_1=A_2\xi^e$, and $A_1A_2=1$.\\
        So, $A_1^2A_2=A_2^3\xi^{e}$, $A_1A_2^2=A_1^3\xi^{2e}$, and $A_0A_1A_2=2.$ Then ${Ind_\gamma}(A_1^2A_2)\;\equiv \;es\; ( \bmod\; 3)$, ${Ind_\gamma}(A_1A_2^2)\;\equiv \;2es\; ( \bmod\; 3)$, and ${Ind_\gamma}(A_0A_1A_2)\;\equiv \;{Ind_\gamma}(2)\; ( \bmod\; 3)$. \\
        Using Theorem (\ref{Suff}) and Theorem (\ref{pbn}), $P(x)=x^r(x^{es}+1)$ is a permutation polynomial over $\mathbbm{F}_q$ if and only if $(r,s)=1, 3 \mid 2^s-1, 3\nmid 2r+es, 2{Ind_\gamma}(A_1^2A_2)\;\not\equiv \;-2r\; ( \bmod\; 3), 2{Ind_\gamma}(A_1A_2^2)\;\not\equiv \;-4r\; ( \bmod\; 3),$ and $2{Ind_\gamma}(A_0A_2^2)\;\not\equiv \;-2r\; ( \bmod\; 3).$
        Now 
        \begin{align*}
            & 2{Ind_\gamma}(A_1^2A_2)\;\equiv \;-2r\; ( \bmod\; 3)\\
            \Leftrightarrow \;&
           es\;\equiv \;-r\; ( \bmod\; 3)\\
            \Leftrightarrow \;&
            3 \mid r+es.
        \end{align*}
        So, $2{Ind_\gamma}(A_1^2A_2)\;\not\equiv \;-2r\; ( \bmod\; 3)$ is equivalent with $3\nmid r+es.$ \\ Similarly, $2{Ind_\gamma}(A_1A_2^2)\;\not\equiv \;-4r\; ( \bmod\; 3)$ is equivalent with $3\nmid r+es$, and $2{Ind_\gamma}(A_0^2A_2)\;\not\equiv \;-2r\; ( \bmod\; 3)$ is equivalent with $3\nmid r+2es.$        
        \end{proof}
      \begin{proposition}\label{prop}
          Let $p$ be an odd prime such that $3\mid p-1$ and $p\nmid 2^{\frac{p-1}{3}}-1$. Then for $0<r<p-1$ and $e \in \mathbbm{N}$, there does not exist any permutation binomial of the form $x^r\{x^{\frac{e(p-1)}{3}}+1\}$ over $\mathbbm{F}_p$ where $(e,3)=1$ and $(r,\frac{p-1}{3})=1$.
      \end{proposition} 
      %\begin{proof}
          %Suppose $p-1=ls$ where $l=3, s=\frac{p-1}{3}$.\\          Let $3\mid Ind_\gamma(2)$ where $\gamma$ is a primitive element of $\mathbbm{F}_p$. Then $ {Ind_\gamma}(2^s)\;\equiv \;0\; ( \bmod\; p-1)$, that is, $2^s=1$ (mod $p$).\\
          %So $p\mid 2^{\frac{p-1}{3}}-1$, which is not possible.\\
          %Hence $3\nmid Ind_\gamma(2)$. From Theorem \ref{PB}, it follows that $x^r\{x^{\frac{e(p-1)}{3}}+1\}$ is not a permutation binomial over $\mathbbm{F}_p$ where $0<r<p-1, $ $(e,3)=1,$ and $(r,\frac{p-1}{3})=1$.
      %\end{proof}
      Proof of Proposition (\ref{prop}) follows from Theorem (\ref{pbn}). By Proposition (\ref{prop}), for $p=7,13,19$; there are no permutation binomials of the form $x^r\{x^{\frac{e(p-1)}{3}}+1\}$ over $\mathbbm{F}_p$, where $(e,3)=1$ and $(r,\frac{p-1}{3})=1$. However, permutation binomials of that form may exist over $\mathbbm{F}_{31}$. 
      \begin{lemma}\label{lem}
          $\xi^{e(j-i)} \frac{A_i}{A_j}=\frac{A_{l-i}}{A_{l-j}}$, for any $i,j$ with $1 \leq i \neq j \leq l-1.$
      \end{lemma}
      \begin{proof}
          For $i=0,1,\cdots,l-1$, We have $A_i=f(\xi^i)=\xi^{ei}+1$ where $\xi$ is a primitive $l$-th root of unity. Then trivially $A_i\neq 0.$\\
          Now $\xi^{e(j-i)} \frac{A_i}{A_j}=\xi^{e(j-i)}\frac{\xi^{ei}+1}{\xi^{ej}+1}=\frac{\xi^{-ei}+1}{\xi^{-ej}+1}=\frac{\xi^{e(l-i)}+1}{\xi^{e(l-j)}+1}=\frac{A_{l-i}}{A_{l-j}}$, for any $i$ and $j$ with $1 \leq i \neq j \leq l-1.$
      \end{proof}
      \begin{theorem}
          Let $p$ be an odd prime, and $q=p^m$ for $m \in \mathbbm{N}$. Assume $l,e,r,s \in \mathbbm{N}$ such that $l(\geq 3)$ is odd, $s$ is even, $(l,e)=1$, $l\mid r+es$, and $q-1=ls$. Then $P(x)=x^r(x^{es}+1)$ is a permutation binomial over $\mathbbm{F}_q$ if and only if $(r,s)=1, p \mid 2^s-1, l\nmid r,$ $\lambda_{l-1}=\{ A_1^s, A_2^s, \cdots , A_{l-1}^s \}$ is a collection of distinct $l$-th root of unity, and ${Ind_\gamma}(A_k)+kr\;\not\equiv \;{Ind_\gamma}(2)\; ( \bmod\; l) \;\forall k=0,1, \cdots, l-1$.
      \end{theorem}
      \begin{proof}
          For the given conditions, from Theorem (\ref{main}) and Theorem (\ref{pbn}), we know that $P(x)$ is a permutation binomial over $\mathbbm{F}_q$ if and only if $(r,s)=1, p \mid 2^s-1, l\nmid 2r+es$ and $ \mu_l=\{ A_0^s, A_1^s\xi^r, \cdots,  A_{l-1}^s\xi^{(l-1)r}\} \text{ is the collection of all distinct } l \text{-th roots of unity}.$\\
          As $l\mid r+es$, then $l\nmid 2r+es$ is equivalent with $l\nmid r$.\\
          For some $i$ and $j$ with $1 \leq i \neq j \leq l-1$, Suppose $A_i^s\xi^{ir} =A_j^s\xi^{jr}$. Then
          \begin{align*}
            & A_i^s\xi^{ir} =A_j^s\xi^{jr} \\
 \Leftrightarrow & \; (A_i/A_j)^s=\xi^{r(j-i)}\\
 \Leftrightarrow & \; \xi^{es(j-i)}(A_i/A_j)^s=\xi^{(r+es)(j-i)}=1 \\
 \Leftrightarrow & \; (A_{l-i}/A_{l-j})^s=1, \text{(using Lemma (\ref{lem}))}\\
 \Leftrightarrow & \; A_{l-i}^s=A_{l-j}^s.
        \end{align*}
          Hence for any $i,j$ with $1 \leq i \neq j \leq l-1$, $A_i^s\xi^{ir} \neq A_j^s\xi^{jr}$ is equivalent with $A_i^s \neq A_j^s$, that is, $\lambda_{l-1}=\{ A_1^s, A_2^s, \cdots , A_{l-1}^s \}$ is a collection of distinct $l$-th root of unity.\\
          For some $k$ with $1 \leq k\leq l-1$, Suppose $A_0^s =A_k^s\xi^{kr}$. Then \begin{align*}
            & A_0^s =A_
            k^s\xi^{kr} \\
 \Leftrightarrow & \; s{Ind_\gamma}(2/A_k)\;\equiv \;krs\; ( \bmod\; q-1) \\
 \Leftrightarrow & \; {Ind_\gamma}(A_k)+kr\;\equiv \;{Ind_\gamma}(2)\; ( \bmod\; l).
        \end{align*}
        Hence for any $k$ with $1 \leq k \leq l-1$, $A_0^s \neq A_k^s\xi^{kr}$ is equivalent with ${Ind_\gamma}(A_k)+kr\;\not\equiv \;{Ind_\gamma}(2)\; ( \bmod\; l)$. Therefore, to show $\mu_l$ is a collection of distinct $l$-th root of unity, it is enough to show that for any $k$ with $1 \leq k \leq l-1$, ${Ind_\gamma}(A_k)+kr\;\not\equiv \;{Ind_\gamma}(2)\; ( \bmod\; l)$ and $\lambda_{l-1}$ is a collection of distinct $l$-th root of unity.
      \end{proof}

    \bibliographystyle{plain}
    \bibliography{bibl.bib}

\begin{thebibliography}{1}

\bibitem{chu2004constructions}
Wensong Chu, Charles~J Colbourn, and Peter Dukes.
\newblock Constructions for permutation codes in powerline communications.
\newblock {\em Designs, Codes and Cryptography}, 32(1):51--64, 2004.

\bibitem{chu2002circular}
Wensong Chu and Solomon~W Golomb.
\newblock Circular tuscan-k arrays from permutation binomials.
\newblock {\em Journal of Combinatorial Theory, Series A}, 97(1):195--202, 2002.

\bibitem{golomb2002periodicity}
Solomon~W Golomb and Oscar Moreno.
\newblock On periodicity properties of costas arrays and a conjecture on permutation polynomials.
\newblock {\em IEEE Transactions on Information Theory}, 42(6):2252--2253, 2002.

\bibitem{levine1987some}
Jack Levine and Richard Chandler.
\newblock Some further cryptographic applications of permutation polynomials.
\newblock {\em Cryptologia}, 11(4):211--218, 1987.

\bibitem{lidl1984permutation}
Rudolf Lidl and Winfried~B M{\"u}ller.
\newblock Permutation polynomials in rsa-cryptosystems.
\newblock In {\em Advances in Cryptology: Proceedings of Crypto 83}, pages 293--301. Springer, 1984.

\bibitem{lidl1994introduction}
Rudolf Lidl and Harald Niederreiter.
\newblock {\em Introduction to finite fields and their applications}.
\newblock Cambridge university press, 1994.

\bibitem{rivest2001permutation}
Ronald~L Rivest.
\newblock Permutation polynomials modulo 2w.
\newblock {\em Finite fields and their applications}, 7(2):287--292, 2001.

\bibitem{sun2005interleavers}
Jing Sun and Oscar~Y Takeshita.
\newblock Interleavers for turbo codes using permutation polynomials over integer rings.
\newblock {\em IEEE Transactions on Information Theory}, 51(1):101--119, 2005.

\bibitem{wang2007cyclotomic}
Qiang Wang.
\newblock Cyclotomic mapping permutation polynomials over finite fields.
\newblock In {\em Sequences, Subsequences, and Consequences: International Workshop, SSC 2007, Los Angeles, CA, USA, May 31-June 2, 2007, Revised Invited Papers}, pages 119--128. Springer, 2007.

\end{thebibliography}
\end{document}